# VARIABLE SELECTION IN SEMIPARAMETRIC REGRESSION MODELING

By Runze Li[1] and Hua Liang[2]

*Pennsylvania State University and University of Rochester*

In this paper, we are concerned with how to select significant variables in semiparametric modeling. Variable selection for semiparametric regression models consists of two components: model selection for nonparametric components and selection of significant variables for the parametric portion. Thus, semiparametric variable selection is much more challenging than parametric variable selection (e.g., linear and generalized linear models) because traditional variable selection procedures including stepwise regression and the best subset selection now require separate model selection for the nonparametric components for each submodel. This leads to a very heavy computational burden. In this paper, we propose a class of variable selection procedures for semiparametric regression models using nonconcave penalized likelihood. We establish the rate of convergence of the resulting estimate. With proper choices of penalty functions and regularization parameters, we show the asymptotic normality of the resulting estimate and further demonstrate that the proposed procedures perform as well as an oracle procedure. A semiparametric generalized likelihood ratio test is proposed to select significant variables in the nonparametric component. We investigate the asymptotic behavior of the proposed test and demonstrate that its limiting null distribution follows a chi-square distribution which is independent of the nuisance parameters. Extensive Monte Carlo simulation studies are conducted to examine the finite sample performance of the proposed variable selection procedures.

**1. Introduction.** Semiparametric regression models retain the virtues of both parametric and nonparametric modeling. Härdle, Liang and Gao [13],

Received September 2005; revised March 2007.

[1]Supported by NSF Grant DMS-03-48869 and partially supported by National Institute on Drug Abuse (NIDA) Grant P50 DA10075.

[2]Supported by NIH/NIAID Grants AI62247 and AI59773.

Ruppert, Wand and Carroll [21] and Yatchew [26] present diverse semiparametric regression models along with their inference procedures and applications. The goal of this paper is to develop effective model selection procedures for a class of semiparametric regression models. Let $Y$ be a response variable and $\{U, \mathbf{X}, \mathbf{Z}\}$ its associated covariates. Further, let $\mu(u, \mathbf{x}, \mathbf{z}) = E(Y|U = u, \mathbf{X} = \mathbf{x}, \mathbf{Z} = \mathbf{z})$. The generalized varying-coefficient partially linear model (GVCPLM) assumes that

$$(1.1) \qquad g\{\mu(u, \mathbf{x}, \mathbf{z})\} = \mathbf{x}^T \boldsymbol{\alpha}(u) + \mathbf{z}^T \boldsymbol{\beta},$$

where $g(\cdot)$ is a known link function, $\boldsymbol{\beta}$ is a vector of unknown regression coefficients and $\boldsymbol{\alpha}(\cdot)$ is a vector consisting of unspecified smooth regression coefficient functions. Model (1.1) is a semiparametric model, $\mathbf{z}^T \boldsymbol{\beta}$ is referred to as the *parametric component* and $\mathbf{x}^T \boldsymbol{\alpha}(u)$ as the *nonparametric component* as $\boldsymbol{\alpha}(\cdot)$ is nonparametric. This semiparametric model retains the flexibility of a nonparametric regression model and has the explanatory power of a generalized linear regression model. Many existing semiparametric or nonparametric regression models are special cases of model (1.1). For instance, partially linear models (see, e.g., Härdle, Liang and Gao [13] and references therein), generalized partially linear models (Severini and Staniswalis [23] and Hunsberger [15]), semivarying-coefficient models (Zhang, Lee and Song [27], Xia, Zhang and Tong [25] and Fan and Huang [8]) and varying coefficient models (Hastie and Tibshirani [14] and Cai, Fan and Li [4]) can be written in the form of (1.1). Thus, the newly proposed procedures provide a general framework of model selection for these existing models.

Variable selection is fundamental in statistical modeling. In practice, a number of variables are available for inclusion in an initial analysis, but many of them may not be significant and should be excluded from the final model in order to increase the accuracy of prediction. Variable selection for the GVCPLM is challenging in that it includes selection of significant variables in the nonparametric component as well as identification of significant variables in the parametric component. Traditional variable selection procedures such as stepwise regression and the best subset variable selection for linear models may be extended to the GVCPLM, but this poses great challenges because for each submodel, it is necessary to choose smoothing parameters for the nonparametric component. This will dramatically increase the computational burden. In an attempt to simultaneously select significant variables and estimate unknown regression coefficients, Fan and Li [9] proposed a family of variable selection procedures for parametric models via nonconcave penalized likelihood. For linear regression models, this family includes bridge regression (Frank and Friedman [12]) and LASSO (Tibshirani [24]). It has been demonstrated that with proper choice of penalty function and regularization parameters, the nonconcave penalized likelihood estimator performs as well as an oracle estimator (Fan and Li [9]). This encourages us to adopt this methodology for semiparametric regression models. In this paper, we propose a class of variable selection procedures for the parametric component of the GVCPLM. We also study the asymptotic properties of the resulting estimator. We illustrate how the rate of convergence of the resulting estimate depends on the regularization parameters and further establish the oracle properties of the resulting estimate. To select significant variables in the nonparametric component of the GVCPLM, we extend generalized likelihood ratio tests (GLRT, Fan, Zhang and Zhang [10]) from fully nonparametric models to semiparametric models. We show the Wilks phenomenon for model (1.1): the limiting null distribution of the proposed GLRT does not depend on the unknown nuisance parameter and it follows a chi-square distribution with diverging degrees of freedom. This allows us to easily obtain critical values for the GLRT using either the asymptotic chi-square distribution or a bootstrap method.

The paper is organized as follows. In Section 2, we first propose a class of variable selection procedures for the parametric component via the nonconcave penalized likelihood approach and then study the sampling properties of the proposed procedures. In Section 3, variable selection procedures are proposed for the nonparametric component using GLRT. The limiting null distribution of the GLRT is derived. Monte Carlo studies and an application involving real data are presented in Section 4. Regularity conditions and technical proofs are presented in Section 5.

**2. Selection of significant variables in the parametric component.** Suppose that $\{U_i, \mathbf{X}_i, \mathbf{Z}_i, Y_i\}$, $i = 1, \ldots, n$, constitute an independent and identically distributed sample and that conditionally on $\{U_i, \mathbf{X}_i, \mathbf{Z}_i\}$, the conditional quasi-likelihood of $Y_i$ is $Q\{\mu(U_i, \mathbf{X}_i, \mathbf{Z}_i), Y_i\}$, where the quasi-likelihood function is defined by

$$Q(\mu, y) = \int_\mu^y \frac{s - y}{V(s)} \, ds$$

for a specific variance function $V(s)$. Throughout this paper, $\mathbf{X}_i$ is $p$-dimensional $\mathbf{Z}_i$ is $d$-dimensional and $U$ is univariate. The methods can be extended for multivariate $U$ in a similar way without any essential difficulty. However, the extension may not be very useful in practice due to the "curse of dimensionality."

2.1. *Penalized likelihood.* Denote by $\ell(\boldsymbol{\alpha}, \boldsymbol{\beta})$ the quasi-likelihood of the collected data $\{(U_i, \mathbf{X}_i, \mathbf{Z}_i, Y_i), i = 1, \ldots, n\}$. That is,

$$\ell(\boldsymbol{\alpha}, \boldsymbol{\beta}) = \sum_{i=1}^n Q[g^{-1}\{\mathbf{X}_i^T \boldsymbol{\alpha}(U_i) + \mathbf{Z}_i^T \boldsymbol{\beta}\}, Y_i].$$

Following Fan and Li [9], define the penalized quasi-likelihood as

$$(2.1) \qquad \mathcal{L}(\boldsymbol{\alpha}, \boldsymbol{\beta}) = \ell(\boldsymbol{\alpha}, \boldsymbol{\beta}) - n \sum_{j=1}^{d} p_{\lambda_j}(|\beta_j|),$$

where $p_{\lambda_j}(\cdot)$ is a prespecified penalty function with a regularization parameter $\lambda_j$, which can be chosen by a data-driven criterion such as cross-validation (CV) or generalized cross-validation (GCV, Craven and Wahba [6]). Note that the penalty functions and regularization parameters are not necessarily the same for all $j$. For example, we wish to keep some important variables in the final model and therefore do not want to penalize their coefficients.

Before we pursue this further, let us briefly discuss how to select the penalty functions. Various penalty functions have been used in the literature on variable selection for linear regression models. Let us take the penalty function to be the $L_0$ penalty, namely, $p_{\lambda_j}(|\beta|) = 0.5\lambda_j^2 I(|\beta| \neq 0)$, where $I(\cdot)$ is the indicator function. Note that $\sum_{j=1}^{d} I(|\beta_j| \neq 0)$ equals the number of nonzero regression coefficients in the model. Hence, many popular variable selection criteria such as AIC (Akaike [1]), BIC (Schwarz [22]) and RIC (Foster and George [11]) can be derived from a penalized least squares problem with the $L_0$ penalty by choosing different values of $\lambda_j$, even though these criteria were motivated by different principles. Since the $L_0$ penalty is discontinuous, it requires an exhaustive search of all possible subsets of predictors to find the solution. This approach is very expensive in computational cost when the dimension $d$ is large. Furthermore, the best subset variable selection suffers from other drawbacks, the most severe of which is its lack of stability, as analyzed, for instance, by Breiman [3].

To avoid the drawbacks of the best subset selection, that is, expensive computational cost and the lack of stability, Tibshirani [24] proposed the LASSO, which can be viewed as the solution of the penalized least squares problem with the $L_1$ penalty, defined by $p_{\lambda_j}(|\beta|) = \lambda_j|\beta|$. Frank and Friedman [12] considered the $L_q$ penalty, $p_{\lambda_j}(|\beta|) = \lambda_j|\beta|^q$, $0 < q < 1$, which yields a "bridge regression." The issue of the selection of the penalty function has been studied in depth by various authors, for instance, Antoniadis and Fan [2]. Fan and Li [9] suggested using the smoothly clipped absolute deviation (SCAD) penalty, defined by

$$p'_{\lambda_j}(\beta) = \lambda_j \left\{ I(\beta \leq \lambda_j) + \frac{(a\lambda_j - \beta)_+}{(a-1)\lambda_j} I(\beta > \lambda_j) \right\}$$

$$\text{for some } a > 2 \text{ and } \beta > 0,$$

with $p_{\lambda_j}(0) = 0$. This penalty function involves two unknown parameters, $\lambda_j$ and $a$. Arguing from a Bayesian statistical point of view, Fan and Li [9] suggested using $a = 3.7$. This value will be used in Section 4.

Since $\boldsymbol{\alpha}(\cdot)$ consists of nonparametric functions, (2.1) is not yet ready for optimization. We must first use local likelihood techniques (Fan and Gijbels [7]) to estimate $\boldsymbol{\alpha}(\cdot)$, then substitute the resulting estimate into (2.1) and finally maximize (2.1) with respect to $\boldsymbol{\beta}$. We can thus obtain a penalized likelihood estimate for $\boldsymbol{\beta}$. With specific choices of penalty function, the resulting estimate of $\boldsymbol{\beta}$ will contain some exact zero coefficients. This is equivalent to excluding the corresponding variables from the final model. We thus achieve the objective of variable selection.

Specifically, we linearly approximate $\alpha_j(v)$ for $v$ in a neighborhood of $u$ by

$$\alpha_j(v) \approx \alpha_j(u) + \alpha_j'(u)(v-u) \equiv a_j + b_j(v-u).$$

Let $\mathbf{a} = (a_1, \ldots, a_p)^T$ and $\mathbf{b} = (b_1, \ldots, b_p)^T$. The local likelihood method is to maximize the local likelihood function

$$(2.2) \qquad \sum_{i=1}^{n} Q[g^{-1}\{\mathbf{a}^T\mathbf{X}_i + \mathbf{b}^T\mathbf{X}_i(U_i-u) + \mathbf{Z}_i^T\boldsymbol{\beta}\}, Y_i]K_h(U_i-u)$$

with respect to $\mathbf{a}$, $\mathbf{b}$ and $\boldsymbol{\beta}$, where $K(\cdot)$ is a kernel function and $K_h(t) = h^{-1}K(t/h)$ is a rescaling of $K$ with bandwidth $h$. Let $\{\tilde{\mathbf{a}}, \tilde{\mathbf{b}}, \tilde{\boldsymbol{\beta}}\}$ be the solution of maximizing (2.2). Then

$$\tilde{\boldsymbol{\alpha}}(u) = \tilde{\mathbf{a}}.$$

As demonstrated in Lemma 2, $\tilde{\boldsymbol{\alpha}}$ is $\sqrt{nh}$-consistent, but its efficiency can be improved by the estimator proposed in Section 3.1. As $\boldsymbol{\beta}$ was estimated locally, the resulting estimate $\tilde{\boldsymbol{\beta}}$ does not have root-$n$ convergence rate. To improve efficiency, $\boldsymbol{\beta}$ should be estimated using global likelihood.

Replacing $\boldsymbol{\alpha}$ in (2.1) by its estimate, we obtain the penalized likelihood

$$(2.3) \qquad \mathcal{L}_P(\boldsymbol{\beta}) = \sum_{i=1}^{n} Q\{g^{-1}(\mathbf{X}_i^T\tilde{\boldsymbol{\alpha}}(U_i) + \mathbf{Z}_i^T\boldsymbol{\beta}), Y_i\} - n\sum_{j=1}^{d} p_{\lambda_j}(|\beta_j|).$$

Maximizing $\mathcal{L}_P(\boldsymbol{\beta})$ results in a penalized likelihood estimator $\hat{\boldsymbol{\beta}}$. The proposed approach is in the same spirit as the one-step backfitting algorithm estimate, although one may further employ the backfitting algorithm method with a full iteration or profile likelihood approach to improve efficiency. The next theorem demonstrates that $\hat{\boldsymbol{\beta}}$ performs as well as an oracle estimator in an asymptotic sense. Compared with fully iterated backfitting algorithms and profile likelihood estimation, the newly proposed method is much less computationally costly and is easily implemented. For high-dimensional $X$- and $Z$-variables, the Hessian matrix of the local likelihood function (2.2) may be nearly singular. To make the resulting estimate stable, one may apply the idea of ridge regression to the local likelihood function. See Cai, Fan and Li [4] for a detailed implementation of ridge regression to local modeling.

2.2. *Sampling properties.* We next study the asymptotic properties of the resulting penalized likelihood estimate. We first introduce some notation. Let $\boldsymbol{\alpha}_0(\cdot)$ and $\boldsymbol{\beta}_0$ denote the true values of $\boldsymbol{\alpha}(\cdot)$ and $\boldsymbol{\beta}$, respectively. Furthermore, let $\boldsymbol{\beta}_0 = (\beta_{10}, \ldots, \beta_{d0})^T = (\boldsymbol{\beta}_{10}^T, \boldsymbol{\beta}_{20}^T)^T$. For ease of presentation and without loss of generality, it is assumed that $\boldsymbol{\beta}_{10}$ consists of all nonzero components of $\boldsymbol{\beta}_0$ and that $\boldsymbol{\beta}_{20} = \mathbf{0}$. Let

$$(2.4) \qquad a_n = \max_{1 \le j \le d} \{|p'_{\lambda_j}(|\beta_{j0}|)|, \beta_{j0} \ne 0\}$$

and

$$b_n = \max_{1 \le j \le d} \{|p''_{\lambda_j}(|\beta_{j0}|)|, \beta_{j0} \ne 0\}.$$

THEOREM 1. *Under the regularity conditions given in Section 5, if $a_n \to 0$, $b_n \to 0$, $nh^4 \to 0$ and $nh^2/\log(1/h) \to \infty$ as $n \to \infty$, then there exists a local maximizer $\hat{\beta}$ of $\mathcal{L}_P(\boldsymbol{\beta})$ defined in (2.3) such that its rate of convergence is $O_P(n^{-1/2} + a_n)$, where $a_n$ is given in (2.4).*

We require further notation in order to present the oracle properties of the resulting penalized likelihood estimate. Define

$$\mathbf{b}_n = \{p'_{\lambda_1}(|\beta_{10}|)\operatorname{sgn}(\beta_{10}), \ldots, p'_{\lambda_s}(|\beta_{s0}|)\operatorname{sgn}(\beta_{s0})\}^T$$

and

$$\Sigma_\lambda = \operatorname{diag}\{p''_{\lambda_1}(|\beta_{10}|), \ldots, p''_{\lambda_s}(|\beta_{s0}|)\},$$

where $s$ is the number of nonzero components of $\boldsymbol{\beta}_0$. Let $\mu_j = \int t^j K(t)\,dt$ and $\nu_j = \int t^j K^2(t)\,dt$ for $j = 0, 1, 2$. Define

$$\rho_l(t) = \frac{\{dg^{-1}(t)/dt\}^l}{\sigma^2 V\{g^{-1}(t)\}} \qquad \text{for } l = 1, 2$$

and let $q_1(x, y) = \rho_1(x)\{y - g^{-1}(x)\}$. Let $R = \boldsymbol{\alpha}_0^T(U)\mathbf{X} + \mathbf{Z}_1^T\boldsymbol{\beta}_{10}$ and

$$(2.5) \qquad \boldsymbol{\Sigma}(u) = E\left[\rho_2(R)\begin{pmatrix} \mathbf{XX}^T & \mathbf{XZ}^T \\ \mathbf{ZX}^T & \mathbf{ZZ}^T \end{pmatrix} \Big| U = u\right].$$

Denote by $\kappa_k$ the $k$th element of $q_1(R, Y)\boldsymbol{\Sigma}^{-1}(u)(\mathbf{X}_1^T, \mathbf{Z}_1^T)^T$ and define

$$\Gamma_1(u) = \sum_{k=1}^{p} \kappa_k E[\rho_2(R)X_k\mathbf{Z}_1 | U = u].$$

THEOREM 2. *Suppose that the regularity conditions given in Section 5 hold and that*

$$\liminf_{n \to \infty} \liminf_{\beta_j \to 0^+} \lambda_{jn}^{-1} p'_{\lambda_{jn}}(|\beta_j|) > 0.$$

*If $\sqrt{n}\lambda_{jn} \to \infty$, $nh^4 \to 0$ and $nh^2/\log(1/h) \to \infty$ as $n \to \infty$, then the root-$n$ consistent estimator $\hat{\boldsymbol{\beta}}$ in Theorem 1 must satisfy $\hat{\boldsymbol{\beta}}_2 = \mathbf{0}$, and $\sqrt{n}(\mathbf{B}_1 + \Sigma_\lambda)\{\hat{\boldsymbol{\beta}}_1 - \boldsymbol{\beta}_{10} + (\mathbf{B}_1 + \Sigma_\lambda)^{-1}\mathbf{b}_n\} \xrightarrow{D} N(0, \Sigma)$, where $\mathbf{B}_1 = [\rho_2(R)\mathbf{Z}_1\mathbf{Z}_1^T]$ and $\Sigma = \mathrm{var}\{q_1(R, Y)\mathbf{Z}_1 - \Gamma_1(U)\}$.*

Theorem 2 indicates that undersmoothing is necessary in order for $\hat{\boldsymbol{\beta}}$ to have root-$n$ consistency and asymptotic normality. This is a standard result in generalized partially linear models; see Carroll et al. [5] for a detailed discussion. Thus, special care is needed for bandwidth selection, as discussed in Section 3.1.

2.3. *Issues arising in practical implementation.* Many penalty functions $p_{\lambda_j}(|\beta_j|)$, including the $L_1$ penalty and the SCAD penalty, are irregular at the origin and may not have a second derivative at some points. Direct implementation of the Newton–Raphson algorithm may be difficult. Following Fan and Li [9], we locally approximate the penalty function by a quadratic function at every step of the iteration, as follows. Given an initial value $\boldsymbol{\beta}^{(0)}$ that is close to the maximizer of the penalized likelihood function, when $\beta_j^{(0)}$ is not very close to 0, the penalty $p_{\lambda_j}(|\beta_j|)$ can be locally approximated by the quadratic function as

$$[p_{\lambda_j}(|\beta_j|)]' = p'_{\lambda_j}(|\beta_j|)\,\mathrm{sgn}(\beta_j) \approx \{p'_{\lambda_j}(|\beta_j^{(0)}|)/|\beta_j^{(0)}|\}\beta_j.$$

Otherwise, set $\hat{\beta}_j = 0$. In other words, for $\beta_j \approx \beta_j^{(0)}$, we have

$$p_{\lambda_j}(|\beta_j|) \approx p_{\lambda_j}(|\beta_j^{(0)}|) + \tfrac{1}{2}\{p'_{\lambda_j}(|\beta_j^{(0)}|)/|\beta_j^{(0)}|\}(\beta_j^2 - \beta_j^{(0)2}).$$

For instance, this local quadratic approximation for the $L_1$ penalty yields

$$|\beta_j| \approx \frac{1}{2}|\beta_j^{(0)}| + \frac{1}{2}\frac{\beta_j^2}{|\beta_j^{(0)}|} \qquad \text{for } \beta_j \approx \beta_j^{(0)}.$$

With the aid of the local quadratic approximation, the Newton–Raphson algorithm can be modified to search for the solution of the penalized likelihood. The convergence of the modified Newton–Raphson algorithm for other statistical settings has been studied by Hunter and Li [16].

*Standard error formula for $\hat{\boldsymbol{\beta}}$.* The standard errors for estimated parameters can be obtained directly because we are estimating parameters and selecting variables at the same time. Following the conventional technique in the likelihood setting, the corresponding sandwich formula can be used as an estimator for the covariance matrix of the estimates $\hat{\boldsymbol{\beta}}$. Specifically, let

$$\ell'(\boldsymbol{\beta}) = \frac{\partial \ell(\tilde{\boldsymbol{\alpha}}, \boldsymbol{\beta})}{\partial \boldsymbol{\beta}}, \qquad \ell''(\boldsymbol{\beta}) = \frac{\partial^2 \ell(\tilde{\boldsymbol{\alpha}}, \boldsymbol{\beta})}{\partial \boldsymbol{\beta}\,\partial \boldsymbol{\beta}^T}$$

and

$$\Sigma_\lambda(\boldsymbol{\beta}) = \mathrm{diag}\left\{\frac{p'_{\lambda_1}(|\beta_1|)}{|\beta_1|}, \ldots, \frac{p'_{\lambda_d}(|\beta_d|)}{|\beta_d|}\right\}.$$

The corresponding sandwich formula is then given by

$$\widehat{\mathrm{cov}}(\hat{\boldsymbol{\beta}}) = \{\ell''(\hat{\boldsymbol{\beta}}) - n\Sigma_\lambda(\hat{\boldsymbol{\beta}})\}^{-1} \widehat{\mathrm{cov}}\{\ell'(\hat{\boldsymbol{\beta}})\}\{\ell''(\hat{\boldsymbol{\beta}}) - n\Sigma_\lambda(\hat{\boldsymbol{\beta}})\}^{-1}.$$

This formula can be shown to be a consistent estimator and will be shown to have good accuracy for moderate sample sizes.

*Choice of $\lambda_j$'s.* We suggest selecting the tuning parameters $\lambda_j$ using data-driven approaches. Similarly to Fan and Li [9], we will employ generalized cross-validation (GCV) to select the $\lambda_j$'s. In the last step of the Newton–Raphson iteration, we may compute the number of effective parameters:

$$e(\lambda_1, \ldots, \lambda_d) = \mathrm{tr}[\{\ell''(\hat{\boldsymbol{\beta}}) - n\Sigma_\lambda(\hat{\boldsymbol{\beta}})\}^{-1}\ell''(\hat{\boldsymbol{\beta}})].$$

The GCV statistic is defined by

$$\mathrm{GCV}(\lambda_1, \ldots, \lambda_d) = \frac{\sum_{i=1}^n D\{Y_i, g^{-1}(\mathbf{X}_i^T\hat{\boldsymbol{\alpha}}(U_i) + \mathbf{Z}_i^T\hat{\boldsymbol{\beta}}(\boldsymbol{\lambda}))\}}{n\{1 - e(\lambda_1, \ldots, \lambda_d)/n\}^2},$$

where $D\{Y, \mu\}$ denotes the deviance of $Y$ corresponding to the model fit with $\boldsymbol{\lambda}$. The minimization problem over a $d$-dimensional space is difficult. However, it is expected that the magnitude of $\lambda_j$ should be proportional to the standard error of the unpenalized maximum pseudo-partial likelihood estimator of $\beta_j$. In practice, we suggest taking $\lambda_j = \lambda SE(\hat{\beta}_j^u)$, where $SE(\hat{\beta}_j^u)$ is the estimated standard error of $\hat{\beta}_j^u$, the unpenalized likelihood estimate. Such a choice of $\lambda_j$ works well from our simulation experience. The minimization problem will thus reduce to a one-dimensional problem and the tuning parameter can be estimated by means of a grid search.

**3. Statistical inferences for nonparametric components.** In this section, we will propose an estimation procedure for $\boldsymbol{\alpha}(\cdot)$ and extend the generalized likelihood ratio test from nonparametric models to model (1.1).

3.1. *Estimation of nonparametric component.* Replacing $\boldsymbol{\beta}$ in (2.2) by its estimate $\hat{\boldsymbol{\beta}}$, we maximize the local likelihood function

$$(3.1) \qquad \sum_{i=1}^n Q[g^{-1}\{\mathbf{a}^T\mathbf{X}_i + \mathbf{b}^T\mathbf{X}_i(U_i - u) + \mathbf{Z}_i^T\hat{\boldsymbol{\beta}}\}, Y_i]K_h(U_i - u)$$

with respect to $\mathbf{a}$ and $\mathbf{b}$. Let $\{\hat{\mathbf{a}}, \hat{\mathbf{b}}\}$ be the solution of maximizing (3.1) and let $\hat{\boldsymbol{\alpha}}(u) = \hat{\mathbf{a}}$. Similarly to Cai, Fan and Li [4], we can show that

$$(nh)^{1/2}\left\{\hat{\boldsymbol{\alpha}}(u) - \boldsymbol{\alpha}_0(u) - \frac{\mu_2}{2}\boldsymbol{\alpha}_0''(u)h^2\right\} \xrightarrow{D} N\left\{\mathbf{0}, \frac{\nu_0}{f(u)}\Sigma_*(u)\right\},$$

where $\Sigma_*(u) = (E[\rho_2\{\boldsymbol{\alpha}_0^T(U)\mathbf{X} + \boldsymbol{\beta}_0^T\mathbf{Z}\}\mathbf{X}\mathbf{X}^T|U=u])^{-1}$ and where $f(u)$ is the density of $U$. Thus, $\hat{\boldsymbol{\alpha}}(u)$ has conditional asymptotic bias $0.5h^2\mu_2\boldsymbol{\alpha}_0''(u) + o_P(h^2)$ and conditional asymptotic covariance $(nh)^{-1}\nu_0\Sigma_*(u)f^{-1}(u) + o_P(\frac{1}{nh})$. From Lemma 2, the asymptotic bias of $\hat{\boldsymbol{\alpha}}$ is the same as that of $\tilde{\boldsymbol{\alpha}}$, while the asymptotic covariance of $\hat{\boldsymbol{\alpha}}$ is smaller than that of $\tilde{\boldsymbol{\alpha}}$.

A theoretic optimal local bandwidth for estimating the elements of $\boldsymbol{\alpha}(\cdot)$ can be obtained by minimizing the conditional mean squared error (MSE) given by

$$E\{\|\hat{\boldsymbol{\alpha}}(u) - \boldsymbol{\alpha}(u)\|^2|\mathbf{Z}, \mathbf{X}\} = \frac{1}{4}h^4\mu_2^2\|\boldsymbol{\alpha}_0''(u)\|^2 + \frac{1}{nh}\frac{\nu_0\operatorname{tr}\{\Sigma_*(u)\}}{f(u)} + o_P\left(h^4 + \frac{1}{nh}\right),$$

where $\|\cdot\|$ is the Euclidean distance. Thus, the ideal choice of local bandwidth is

$$\hat{h}_{\mathrm{opt}} = \left\{\frac{\nu_0\operatorname{tr}\{\Sigma_*(u)\}}{f(u)\mu_2^2\|\boldsymbol{\alpha}_0''(u)\|^2}\right\}^{1/5} n^{-1/5}.$$

With expressions for the asymptotic bias and variance, we can also derive a theoretic or data-driven global bandwidth selector by utilizing the existing bandwidth selection techniques for the canonical univariate nonparametric model, such as the substitution method (see, e.g., Ruppert, Sheather and Wand [20]). For the sake of brevity, we omit the details here.

As usual, the optimal bandwidth will be of order $n^{-1/5}$. This does not satisfy the condition in Theorems 1 and 2. A good bandwidth is generally generated by $\hat{h}_{\mathrm{opt}} \times n^{-2/15} = O(n^{-1/3})$.

In order for the resulting variable selection procedures to possess an oracle property, the bandwidth must satisfy the conditions $nh^4 \to 0$ and $nh^2/(\log n)^2 \to \infty$. The aforementioned order of bandwidth satisfies these requirements. This enables us to easily choose a bandwidth either by data-driven procedures or by an asymptotic theory-based method.

3.2. *Variable selection for the nonparametric component.* After obtaining nonparametric estimates for $\{\alpha_1(\cdot), \ldots, \alpha_p(\cdot)\}$, it is of interest to select significant $X$-variables. For linear regression models, one conducts an $F$-test at each step of the traditional backward elimination or forward addition procedures. The purpose of variable selection may be achieved by a sequence of $F$-tests. Following the strategy of the traditional variable selection procedure, one may apply the backward elimination procedure to select significant $x$-variables. In each step of the backward elimination procedure, we essentially test the following hypothesis

$$H_0 : \alpha_{j_1}(u) = \cdots = \alpha_{j_k}(u) = 0 \quad \text{versus} \quad H_1 : \text{not all } \alpha_{j_l}(u) \neq 0$$

for some $\{j_1, \ldots, j_k\}$, a subset of $\{1, \ldots, p\}$. The purpose of variable selection may be achieved by a sequence of such tests. For ease of presentation, we here consider the following hypothesis:

$$(3.2) \quad H_0 : \alpha_1(u) = \cdots = \alpha_p(u) = 0 \quad \text{versus} \quad H_1 : \text{not all } \alpha_j(u) \neq 0.$$

The proposed idea is also applicable to more general cases.

Let $\widehat{\boldsymbol{\alpha}}(u)$ and $\widehat{\boldsymbol{\beta}}$ be the estimators of $\boldsymbol{\alpha}(u)$ and $\boldsymbol{\beta}$ under the alternative hypothesis, respectively, and let $\bar{\boldsymbol{\beta}}$ be the estimators of $\boldsymbol{\beta}$ under the null hypothesis. Define

$$\mathcal{R}(H_1) = \sum_{i=1}^{n} Q\{g^{-1}(\widehat{\boldsymbol{\alpha}}^T(U_i)\mathbf{X}_i^T + \mathbf{Z}_i^T\widehat{\boldsymbol{\beta}}), Y_i\}$$

and

$$\mathcal{R}(H_0) = \sum_{i=1}^{n} Q\{g^{-1}(\mathbf{Z}_i^T\bar{\boldsymbol{\beta}}), Y_i\}.$$

Following Fan, Zhang and Zhang [10], we define a generalized quasi-likelihood ratio test (GLRT) statistic

$$T_{\text{GLR}} = r_K\{\mathcal{R}(H_1) - \mathcal{R}(H_0)\},$$

where

$$r_K = \left\{K(0) - 0.5 \int K^2(u)\, du\right\} \left\{\int \{K(u) - 0.5K * K(u)\}\, du\right\}^{-1}.$$

THEOREM 3. *Suppose that the regularity conditions given in Section 5 hold and that $nh^8 \to 0$ and $nh^2/(\log n)^2 \to \infty$. Under $H_0$ in (3.2), the test statistic $T_{\text{GLR}}$ has an asymptotic $\chi^2$ distribution with $df_n$ degrees of freedom, in the sense of Fan, Zhang and Zhang [10], where $df_n = r_K p|\Omega|\{K(0) - 0.5 \int K^2(u)\, du\}/h$ and where $|\Omega|$ denotes the length of the support of $U$.*

Theorem 3 reveals a new Wilks phenomenon for semiparametric inference and extends the generalized likelihood ratio theory (Fan, Zhang and Zhang [10]) for semiparametric modeling. We will also provide empirical justification for the null distribution. Similarly to Cai, Fan and Li [4], the null distribution of $T_{\text{GLR}}$ can be estimated using Monte Carlo simulation or a bootstrap procedure. This usually provides a better estimate than the asymptotic null distribution since the degrees of freedom tend to infinity and the results in Fan, Zhang and Zhang [10] only give the main order of the degrees of freedom.

**4. Simulation study and application.** In this section, we conduct extensive Monte Carlo simulations in order to examine the finite sample performance of the proposed procedures.

The performance of the estimator $\hat{\boldsymbol{\alpha}}(\cdot)$ will be assessed by using the square root of average square errors (RASE)

$$\text{(4.1)} \qquad \text{RASE} = \left\{ n_{\text{grid}}^{-1} \sum_{k=1}^{n_{\text{grid}}} \|\hat{\boldsymbol{\alpha}}(u_k) - \boldsymbol{\alpha}(u_k)\|^2 \right\}^{1/2},$$

where $\{u_k, k = 1, \ldots, n_{\text{grid}}\}$ are the grid points at which the functions $\{\hat{\alpha}_j(\cdot)\}$ are evaluated. In our simulation, the Epanechnikov kernel $K(u) = 0.75(1 - u^2)_+$ and $n_{\text{grid}} = 200$ are used.

In an earlier version of this paper (Li and Liang [17]), we assessed the performance of the proposed estimation procedure for $\boldsymbol{\beta}$ without the task of variable selection and concluded that the proposed estimation procedures performs well. We have since further tested the accuracy of the proposed standard error formula and found that it works fairly well. In this section, we focus on the performance of the proposed variable selection procedures. The prediction error is defined as the average error in the prediction of the dependent variable given the independent variables for future cases that are not used in the construction of a prediction equation. Let $\{U^*, \mathbf{X}^*, \mathbf{Z}^*, Y^*\}$ be a new observation from the GVCPLM model (1.1). The prediction error for model (1.1) is then given by

$$\text{PE}(\hat{\boldsymbol{\alpha}}, \hat{\boldsymbol{\beta}}) = E\{Y^* - \hat{\mu}(U^*, \mathbf{X}^*, \mathbf{Z}^*)\}^2,$$

where the expectation is a conditional expectation given the data used in constructing the prediction procedure. The prediction error can be decomposed as

$$\text{PE}(\hat{\boldsymbol{\alpha}}, \hat{\boldsymbol{\beta}}) = E\{Y^* - \mu(U^*, \mathbf{X}^*, \mathbf{Z}^*)\}^2 + E\{\hat{\mu}(U^*, \mathbf{X}^*, \mathbf{Z}^*) - \mu(U^*, \mathbf{X}^*, \mathbf{Z}^*)\}^2.$$

The first component is the inherent prediction error due to noise. The second component is due to the lack of fit with an underlying model. This component is called the *model error*. Note that $\hat{\boldsymbol{\alpha}}, \hat{\boldsymbol{\beta}}$ provide a consistent estimate and that $\mu(U^*, \mathbf{X}^*, \mathbf{Z}^*) = g^{-1}\{\mathbf{Z}^{*T}\boldsymbol{\alpha}(U^*) + \mathbf{Z}^{*T}\boldsymbol{\beta}\}$. By means of a Taylor expansion, we have the approximation

$$\hat{\mu}(U^*, \mathbf{X}^*, \mathbf{Z}^*)$$
$$\approx \mu(U^*, \mathbf{X}^*, \mathbf{Z}^*) + \dot{g}^{-1}\{\mathbf{X}^{*T}\boldsymbol{\alpha}(U^*) + \mathbf{Z}^{*T}\boldsymbol{\beta}\}\mathbf{X}^{*T}\{\hat{\boldsymbol{\alpha}}(U^*) - \boldsymbol{\alpha}(U^*)\}$$
$$+ \dot{g}^{-1}\{\mathbf{X}^{*T}\boldsymbol{\alpha}(U^*) + \mathbf{Z}^{*T}\boldsymbol{\beta}\}\mathbf{Z}^{*T}(\hat{\boldsymbol{\beta}} - \boldsymbol{\beta}),$$

where $\dot{g}^{-1}(t) = dg^{-1}(t)/dt$. Therefore, the model error can be approximated by

$$E[\dot{g}^{-1}\{\mathbf{X}^{*T}\boldsymbol{\alpha}(U^*) + \mathbf{Z}^{*T}\boldsymbol{\beta}\}]^2([\mathbf{X}^{*T}\{\hat{\boldsymbol{\alpha}}(U^*) - \boldsymbol{\alpha}(U^*)\}]^2 + [\mathbf{Z}^{*T}(\hat{\boldsymbol{\beta}} - \boldsymbol{\beta})]^2$$
$$+ [\mathbf{X}^{*T}\{\hat{\boldsymbol{\alpha}}(U^*) - \boldsymbol{\alpha}(U^*)\}] \times [\mathbf{Z}^{*T}(\hat{\boldsymbol{\beta}} - \boldsymbol{\beta})]).$$

The first component is the inherent model error due to the lack of fit of the nonparametric component $\boldsymbol{\alpha}_0(t)$, the second is due to the lack of fit of the parametric component and the third is the cross product between the first two components. Thus, we define generalized mean square error (GMSE) for the parametric component as

$$(4.2) \qquad \text{GMSE}(\hat{\boldsymbol{\beta}}) = E[\mathbf{Z}^{*T}(\hat{\boldsymbol{\beta}} - \boldsymbol{\beta})]^2 = (\hat{\boldsymbol{\beta}} - \boldsymbol{\beta}) E(\mathbf{Z}^* \mathbf{Z}^{*T})(\hat{\boldsymbol{\beta}} - \boldsymbol{\beta})$$

and use the GMSE to assess the performance of the newly proposed variable selection procedures for the parametric component.

EXAMPLE 4.1. In this example, we consider a semivarying Poisson regression model. Given $(U, \mathbf{X}, \mathbf{Z})$, $Y$ has a Poisson distribution with mean function $\mu(U, \mathbf{X}, \mathbf{Z})$ where

$$\mu(U, \mathbf{X}, \mathbf{Z}) = \exp\{\mathbf{X}^T \boldsymbol{\alpha}(U) + \mathbf{Z}^T \boldsymbol{\beta}\}.$$

In our simulation, we take $U \sim U(0,1)$, $\mathbf{X} = (X_1, X_2)^T$ with $X_1 \equiv 1$ and $X_2 \sim N(0,1)$, $\alpha_1(u) = 5.5 + 0.1 \exp(2u - 1)$ and $\alpha_2(u) = 0.8u(1 - u)$. Furthermore, $\boldsymbol{\beta} = [0.3, 0.15, 0, 0, 0.2, 0, 0, 0, 0, 0]^T$ and $\mathbf{Z}$ has a 10-dimensional normal distribution with zero mean and covariance matrix $(\sigma_{ij})_{10 \times 10}$ with $\sigma_{ij} = 0.5^{|i-j|}$. In our simulation, we take the sample size $n = 200$ and bandwidth $h = 0.125$.

*Performance of procedures for $\boldsymbol{\beta}$.* Here, we compare the variable selection procedures for $\mathbf{Z}$. One may generalize the traditional subset selection criteria for linear regression models to the GVCPLM by taking the penalty function in (2.1) to be the $L_0$ penalty. Specifically, $p_{\lambda_j}(|\beta|) = 0.5\lambda_j^2 I(|\beta| \neq 0)$. We will refer to the AIC, BIC and RIC as the penalized likelihood with the $L_0$ penalty with $\lambda_j = \sqrt{2/n}, \sqrt{\log(n)/n}$ and $\sqrt{2\log(d)/n}$, respectively.

TABLE 1
*Comparisons of variable selection*

| Penalty | Poisson with $h = 0.125$ | | | Logistic with $h = 0.3$ | | |
| | RGMSE Median(MAD) | C | I | RGMSE Median(MAD) | C | I |
|---|---|---|---|---|---|---|
| SCAD | 0.3253 (0.2429) | 6.8350 | 0 | 0.5482 (0.3279) | 6.7175 | 0 |
| $L_1$ | 0.8324 (0.1651) | 4.9650 | 0 | 0.7247 (0.2024) | 5.3625 | 0 |
| AIC | 0.7118 (0.2228) | 5.6825 | 0 | 0.8353 (0.1467) | 5.7225 | 0 |
| BIC | 0.3793 (0.2878) | 6.7400 | 0 | 0.5852 (0.3146) | 6.9100 | 0 |
| RIC | 0.4297 (0.2898) | 6.6475 | 0 | 0.6665 (0.2719) | 6.7100 | 0 |
| Oracle | 0.2750 (0.1983) | 7 | 0 | 0.5395 (0.3300) | 7 | 0 |

TABLE 2
*Computing times*

|  | Penalty | $d = 8$ | $d = 9$ | $d = 10$ |
|---|---|---|---|---|
| Poisson | SCAD | 0.0485 (0.0145) | 0.0584 (0.0155) | 0.0471 (0.0132) |
|  | $L_1$ | 0.0613 (0.0145) | 0.0720 (0.0217) | 0.0694 (0.0188) |
|  | BIC | 0.8255 (0.1340) | 2.1558 (0.6836) | 4.6433 (1.3448) |
| Logistic | SCAD | 2.7709 (0.1606) | 2.8166 (0.1595) | 2.8337 (0.1449) |
|  | $L_1$ | 8.1546 (0.8931) | 7.9843 (0.9196) | 8.1952 (0.9491) |
|  | BIC | 61.5723 (1.4404) | 131.8402 (2.6790) | 280.0237 (6.6325) |

Since the $L_0$ penalty is discontinuous, we search over all possible subsets to find the correspond solutions. Thus, these procedures will be referred to as best subset variable selection. We compare the performance of the penalized likelihood with the $L_1$ penalty and the SCAD penalty with the best subset variable selection in terms of GMSE and model complexity. Let us define *relative GMSE* to be the ratio of GMSE of a selected final model to that of the full model. The median of relative GMSE over the 400 simulations, along with the median of absolute deviation divided by a factor of 0.6745, is displayed in the column of Table 1 labeled "RGMSE." The average number of 0 coefficients is also given in Table 1, where the column labeled "C" gives the average number of the seven true zero coefficients, correctly set to zero and the column labeled "I" gives the average number of the three true nonzeros incorrectly set to zero. In Table 1, "Oracle" stands for the oracle estimate computed by using the true model $g\{E(y|u, \mathbf{x}, \mathbf{z})\} = \mathbf{x}^T \boldsymbol{\alpha}(u) + \beta_1 z_1 + \beta_2 z_2 + \beta_5 z_5$. According to Table 1, the performance of the SCAD is close to that of the oracle procedure in terms of model error and model complexity, and it performs better than penalized likelihood with the $L_1$ and best subset variable selection using the AIC and RIC. The performance of the SCAD is similar to that of best subset variable selection with BIC. However, best subset variable selection demands much more computation. To illustrate this, we compare the computing time for each procedure. Table 2 includes the average and standard deviation of computing times over 50 Monte Carlo simulations for $d = 8$, 9 and 10. For $d = 8$ and 9, $\beta_1 = 0.3$, $\beta_2 = 0.15$, $\beta_5 = 0.2$ and other components of $\boldsymbol{\beta} = 0$; $\mathbf{Z}$ has multivariate normal distribution with zero mean and the same covariance structure as that for $d = 10$; $U$, $\mathbf{X}$ and $\boldsymbol{\alpha}(U)$ are the same as those for $d = 10$. We include only the computing time for BIC in Table 2. Computing time for AIC and RIC is almost identical to that for BIC. It is clear from Table 2 that BIC needs much more computing time than the SCAD and $L_1$ and that it exponentially increases as $d$ increases.

*Performance of procedures for* $\boldsymbol{\alpha}(u)$. It is of interest to assess the impact of estimation of $\boldsymbol{\beta}$ on the estimation of $\boldsymbol{\alpha}(\cdot)$. To this end, we consider two scenarios: one is to estimate $\boldsymbol{\alpha}(\cdot)$ using the proposed backfitting algorithm, and the other is to estimate $\boldsymbol{\alpha}(\cdot)$ with the true value of $\boldsymbol{\beta}$. The plot of one RASE versus the other is depicted in Figure 1(a), from which it can be seen that the estimate $\hat{\boldsymbol{\alpha}}$ using the proposed backfitting algorithm performs as well as if we knew the true value of $\boldsymbol{\beta}$. This is consistent with our theoretic analysis because $\hat{\boldsymbol{\beta}}$ is root-$n$ consistent and this convergence rate is faster than the convergence rate of a nonparametric estimate.

We now assess the performance of the test procedures proposed in Section 3. Here, we consider the null hypothesis

$$H_0 : \alpha_2(u) = 0 \quad \text{versus} \quad H_1 : \alpha_2(u) \neq 0.$$

We first examine whether the finite sample null distribution of the proposed GLRT is close to a chi-square distribution. To this end, we conduct 1000 bootstraps to obtain the null distribution of the proposed GLRT. The kernel density estimate of the null distribution is depicted in Figure 1(c), in which the solid line corresponds to the estimated density function and the dotted line to a density of the $\chi^2$-distribution with degrees of freedom approximately equaling the sample mean of the bootstrap sample. From Figure 1(c), the finite sample null distribution is quite close to a chi-square distribution.

We next examine the Type I error rate and power of the proposed GLRT. The power functions are evaluated under a sequence of alternative models indexed by $\delta$:

$$H_1 : \alpha_2(u) = \delta \times 0.8u(1-u).$$

Figure 1(e) depicts four power functions based on 400 simulations at four different significance levels: 0.25, 0.1, 0.05 and 0.01. When $\delta = 0$, the special alternative collapses into the null hypothesis. The powers at $\delta = 0$ for these four significance levels are 0.2250, 0.0875, 0.05 and 0.0125. This shows that the bootstrap method gives a correct Type I error rates. The power functions increase rapidly as $\delta$ increases. This shows that the proposed GLRT works well.

EXAMPLE 4.2. In this example, we consider a semivarying coefficient logistic regression model. Given $(U, \mathbf{X}, \mathbf{Z})$, $Y$ has a Bernoulli distribution with success probability $p(U, \mathbf{X}, \mathbf{Z})$, where

$$p(U, \mathbf{X}, \mathbf{Z}) = \exp\{\mathbf{X}^T \boldsymbol{\alpha}(U) + \mathbf{Z}^T \boldsymbol{\beta}\} / [1 + \exp\{\mathbf{X}^T \boldsymbol{\alpha}(U) + \mathbf{Z}^T \boldsymbol{\beta}\}].$$

In our simulation, $U, \mathbf{X}, \mathbf{Z}$ are the same as those in Example 4.1, but the coefficient functions are taken to be

$$\alpha_1(u) = \exp(2u - 1), \qquad \alpha_2(u) = 2\sin^2(2\pi u)$$

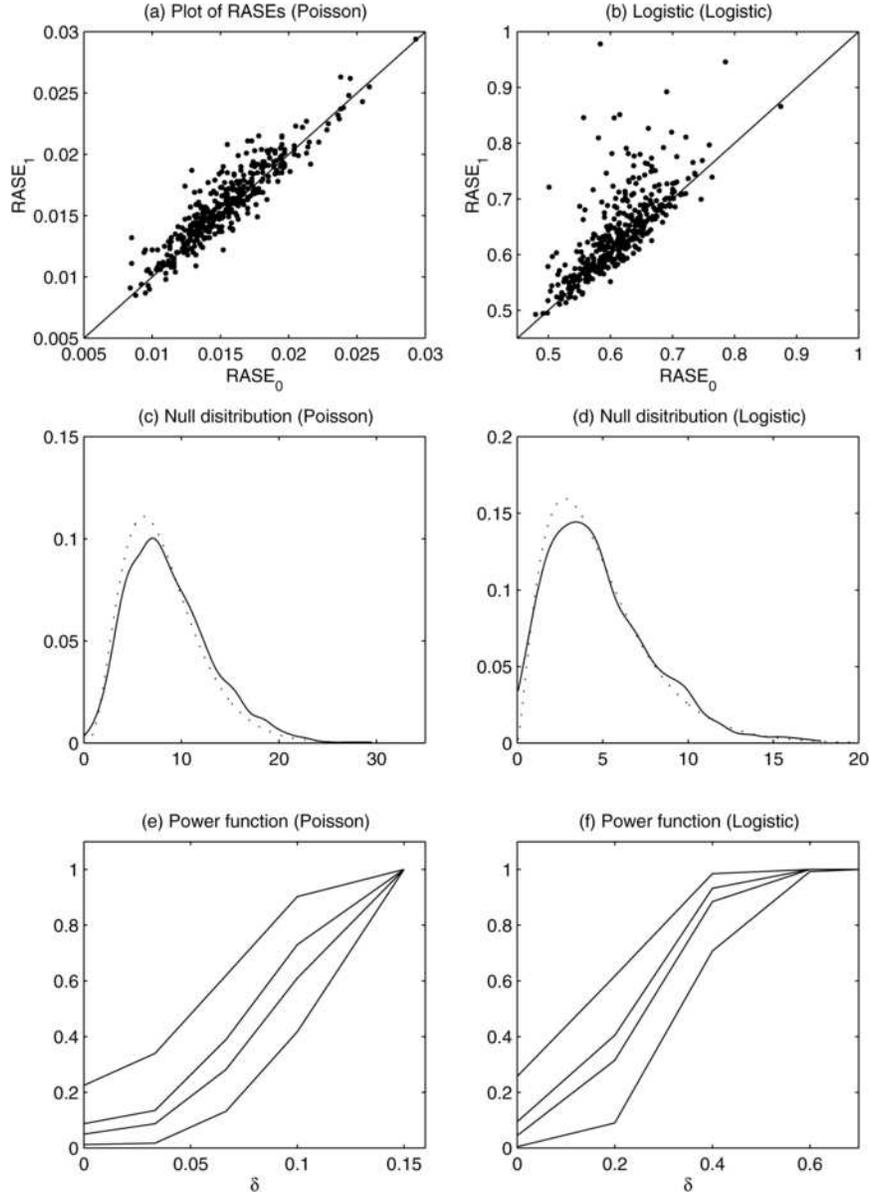

FIG. 1. *Plots for Examples 4.1 (left panel) and 4.2 (right panel). In* (a) *and* (b), RASE$_0$ *stands for the RASE of* $\hat{\boldsymbol{\alpha}}(u)$ *with the true* $\boldsymbol{\beta}$ *and* RASE$_1$ *stands for the RASE of* $\hat{\boldsymbol{\alpha}}(u)$ *using the backfitting algorithm. In* (c) *and* (d), *the solid lines correspond to the estimated null density and the dotted lines to the density of the* $\chi^2$-*distribution with* $df_n$ *being the mean of the bootstrap sample.* (e) *and* (f) *are power functions of the GLRT at levels 0.25, 0.10, 0.05 and 0.01.*

and $\boldsymbol{\beta} = [3, 1.5, 0, 0, 2, 0, 0, 0, 0, 0]^T$. We conduct 400 simulations and in each simulation, the sample size $n$ is set at 1000 and the bandwidth at $h = 0.3$.

*Performance of procedures for $\boldsymbol{\beta}$.* We investigate the performance of the proposed variable selection procedures. Simulation results are summarized in the rightmost column of Table 1, from which we can see that the SCAD performs the best and that its performance is very close to that of the oracle procedure. We employ the same strategy as in Example 4.1 to compare computing times of each variable selection procedure. The mean and standard deviation of computing time are given in the bottom row of Table 2, from which it can be seen that the computing time for the best subset variable selection procedure increases exponentially as the dimension of $\boldsymbol{\beta}$ increases, while this is not the case for penalized likelihood with the SCAD penalty and the $L_1$ penalty.

*Performance of procedures for $\boldsymbol{\alpha}(u)$.* We employ RASE to assess the performance of $\hat{\boldsymbol{\alpha}}(u)$. Figure 1(b) plots the RASE of $\hat{\boldsymbol{\alpha}}(\cdot)$ using the proposed backfitting algorithm against that of $\hat{\boldsymbol{\alpha}}(\cdot)$ using the true value of $\boldsymbol{\beta}$. The performance of the backfitting algorithm is quite close to that using the true value of $\boldsymbol{\beta}$.

We next examine the performance of the proposed GLRT for logistic regression. Here, we consider the null hypothesis

$$H_0 : \alpha_2(u) = 0 \quad \text{versus} \quad H_1 : \alpha_2(u) \neq 0.$$

The estimated density of null distribution is depicted in Figure 1(d), from which we can see that it is close to a $\chi^2$-distribution. The power functions are evaluated under a sequence of alternative models indexed by $\delta$:

$$H_1 : \alpha_2(u) = \delta \times 2 \sin^2(2\pi u).$$

The power functions are depicted in Figure 1(d), from which it can be seen that the power functions increase rapidly as $\delta$ increases.

EXAMPLE 4.3. We now apply the methodology proposed in this paper to the analysis of a data set compiled by the General Hospital Burn Center at the University of Southern California. The binary response variable $Y$ is 1 for those victims who survived their burns and 0 otherwise, the variable $U$ in this application represents *age* and fourteen other covariates were considered. We first employ a generalized varying-coefficient model (Cai, Fan and Li [4]) to fit the data by allowing the coefficients of all fourteen covariates to be age-dependent. Based on the resulting estimates and standard errors, we consider a generalized varying-coefficient partially linear model for binary response,

$$(4.3) \qquad \text{logit}\{E(Y | U = u, \mathbf{X} = \mathbf{x}, \mathbf{Z} = \mathbf{z})\} = \mathbf{x}^T \boldsymbol{\alpha}(u) + \mathbf{z}^T \boldsymbol{\beta},$$

where $X_1 \equiv 1$ and $\alpha_1(u)$ is the intercept function. Other covariates are as follows. $X_2$ stands for log(burn area + 1), $X_3$ for prior respiratory disease (coded by 0 for none and 1 for yes), $Z_1$ for gender (coded by 0 for male and 1 for female), $Z_2$ for days injured prior to admission date (coded by 0 for one or more days and 1 otherwise), $Z_3$ for airway edema (coded by 0 for not present and 1 for present), $Z_4$ for sootiness (coded by 1 for yes and 0 for no), $Z_5$ for partial pressure of oxygen, $Z_6$ for partial pressure of carbon dioxide, $Z_7$ for pH (acidity) reading, $Z_8$ for percentage of CbHg, $Z_9$ for oxygen supply (coded by 0 for normal and 1 for abnormal), $Z_{10}$ for carbon dioxide status (coded by 0 for normal and 1 for abnormal), $Z_{11}$ for acid status coded by (0 for normal and 1 for abnormal) and $Z_{12}$ for hemo status (coded by 0 for normal and 1 for abnormal).

In this demonstration, we are interested in studying how the included covariates affect survival probabilities for victims in different age groups. We first employ a multifold cross-validation method to select a bandwidth. We partition the data into $K$ groups. For each $j$, $k = 1, \dots, K$, we fit the data to model (4.3), excluding data in the $k$th group, denoted by $\mathcal{D}_k$. The deviance is computed. This leads to a cross-validation criterion,

$$CV(h) = \sum_{k=1}^{K} \sum_{i \in \mathcal{D}_k} D\{y_i, \hat{\mu}_{-k}(u_i, \mathbf{x}_i, \mathbf{z}_i)\},$$

where $D(y, \hat{\mu})$ is the deviance of the Bernoulli distribution, $\hat{\mu}_{-k}(u_i, \mathbf{x}_i, \mathbf{z}_i)$ is the fitted value of $Y_i$, that is, $\text{logit}^{-1}\{\mathbf{x}_i^T \hat{\boldsymbol{\alpha}}_{-k}(u_i) + \mathbf{z}_i^T \hat{\beta}_{-k}\}$, and $\hat{\boldsymbol{\alpha}}_{-k}(\cdot)$ and $\hat{\beta}_{-k}$ are estimated without including data from $\mathcal{D}_k$. In our implementation, we set $K = 10$. Figure 2(a) depicts the plot of cross-validation scores over the bandwidth. The selected bandwidth is 48.4437. With the selected bandwidth, the resulting estimate of $\boldsymbol{\alpha}(u)$ is depicted in Figures 2(b), (c) and (d). From the plot of $\hat{\alpha}_3(u)$ in Figure 2(d), we see that the 95% pointwise confidence interval almost covers zero. Thus, it is of interest to test whether or not $X_3$ is significant. To this end, we employ the semiparametric generalized likelihood ratio test procedure for the following hypothesis:

$$H_0 : \alpha_3(u) = 0 \quad \text{versus} \quad H_1 : \alpha_3(u) \neq 0.$$

The resulting generalized likelihood ratio test for this problem is 15.7019, with a P value of 0.015, based on 1,000 bootstrap samples. Thus, the covariate *prior respiratory disease* is significant at level 0.05. The result also implies that the generalized likelihood ratio test is quite powerful as the resulting estimate of $\alpha_3(u)$ only slightly deviates from 0.

We next select significant $z$-variables. We apply the SCAD procedure proposed in Section 2 to the data. The tuning parameter $\lambda$ is chosen by minimizing the GCV scores. The selected $\lambda$ equals 0.4226. With this selected tuning parameter, the SCAD procedure yields a model with only three $z$-variables: $Z_3$, $Z_5$ and $Z_7$. Their estimates and standard errors are $-1.9388(0.4603)$, $-0.0035(0.0054)$ and $-0.0007(0.0006)$, respectively. As a result, we recommend the following model:

$$\hat{Y} = \hat{\alpha_1}(U) + \hat{\alpha_2}(U)X_2 + \hat{\alpha_3}(U)X_3 - 1.9388Z_3 - 0.0035Z_5 - 0.0007Z_7,$$

where the $\hat{\alpha}(U)$'s and their 95% confidence intervals are plotted in Figure 2.

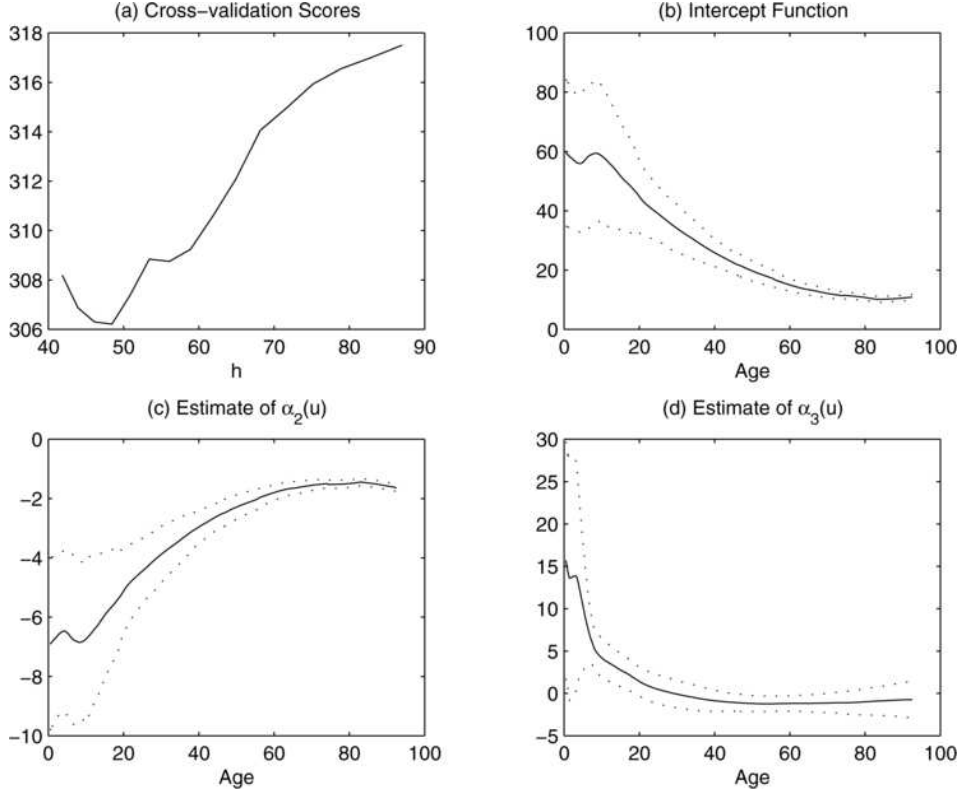

FIG. 2.  *Plots for Example 4.3.*

**5. Proofs.**  For simplicity of notation, in this section, we absorb $\sigma^2$ into $V(\cdot)$, so that the variance of $Y$ given $(U, \mathbf{X}, \mathbf{Z})$ is $V\{\mu(U, \mathbf{X}, \mathbf{Z})\}$. Define $q_\ell(x, y) = (\partial^\ell/\partial x^\ell)Q\{g^{-1}(x), y\}$ for $\ell = 1, 2, 3$. Then

$$(5.1) \qquad \begin{aligned} q_1(x, y) &= \{y - g^{-1}(x)\}\rho_1(x), \\ q_2(x, y) &= \{y - g^{-1}(x)\}\rho_1'(x) - \rho_2(x), \end{aligned}$$

where $\rho_\ell(t) = \{\frac{dg^{-1}(t)}{dt}\}^\ell / V\{g^{-1}(t)\}$ was introduced in Section 2. In the following regularity conditions, $u$ is a generic argument for Theorem 2 and the condition must hold *uniformly* in $u$ for Theorems 1–3.

*Regularity conditions*:

(i) The function $q_2(x, y) < 0$ for $x \in \mathbb{R}$ and $y$ in the range of the response variable.

(ii) The random variable $U$ has bounded support $\Omega$. The elements of the function $\boldsymbol{\alpha}_0''(\cdot)$ are continuous in $u \in \Omega$.

(iii) The density function $f(u)$ of $U$ has a continuous second derivative.

(iv) The functions $V''(\cdot)$ and $g'''(\cdot)$ are continuous.

(v) With $R = \boldsymbol{\alpha}_0^T(U)\mathbf{X} + \mathbf{Z}^T\boldsymbol{\beta}_0$, $E\{q_1^2(R, Y)|U = u\}$, $E\{q_1^2(R, Y)\mathbf{Z}|U = u\}$ and $E\{q_1^2(R, Y)\mathbf{Z}\mathbf{Z}^T|U = u\}$ are twice differentiable in $u \in \Omega$. Moreover, $E\{q_2^2(R, Y)\} < \infty$ and $E\{q_1^{2+\delta}(R, Y)\} < \infty$ for some $\delta > 2$.

(vi) The kernel $K$ is a symmetric density function with bounded support.

(vii) The random vector $\mathbf{Z}$ is assumed to have bounded support.

Condition (i) is imposed so that the local likelihood is concave in the parameters, which ensures the uniqueness of the solution. Conditions (vi) and (vii) are imposed just to simplify the proofs; they can be weakened significantly at the expense of lengthy proofs. In our proofs, we will repeatedly use the following lemma, a direct result of Mack and Silverman [18].

LEMMA 1. *Let $(\mathbf{X}_1, Y_1), \ldots, (\mathbf{X}_n, Y_n)$ be i.i.d. random vectors, where the $Y_i$'s are scalar random variables. Assume further that $E|Y|^r < \infty$ and that $\sup_{\mathbf{x}} \int |y|^r f(\mathbf{x}, y)\,dy < \infty$, where $f$ denotes the joint density of $(\mathbf{X}, Y)$. Let $K$ be a bounded positive function with bounded support, satisfying a Lipschitz condition. Then*

$$\sup_{\mathbf{x} \in D} \left| n^{-1} \sum_{i=1}^{n} \{K_h(\mathbf{X}_i - \mathbf{x})Y_i - E[K_h(\mathbf{X}_i - \mathbf{x})Y_i]\} \right| = O_P\left[\left\{\frac{nh}{\log(1/h)}\right\}^{-1/2}\right],$$

*provided that $n^{2\varepsilon-1}h \to \infty$ for some $\varepsilon < 1 - r^{-1}$.*

To establish asymptotic properties of $\hat{\boldsymbol{\beta}}$, we first study the asymptotic behaviors of $\tilde{\mathbf{a}}, \tilde{\mathbf{b}}$ and $\tilde{\boldsymbol{\beta}}$. Let us introduce some notation. Let $\bar{\alpha}_i = \bar{\alpha}_i(u) = \mathbf{X}^T\boldsymbol{\alpha}_0(u) + \mathbf{Z}_i^T\boldsymbol{\beta}_0 + (U_i - u)\mathbf{X}_i^T\boldsymbol{\alpha}_0'(u)$. Write $\mathbf{X}_i^* = (\mathbf{X}_i^T, (U_i - u)\mathbf{X}_i^T/h, \mathbf{Z}_i^T)^T$,

$$A(\mathbf{X}, \mathbf{Z}) = \begin{pmatrix} \mathbf{X}\mathbf{X}^T & \mathbf{0}^T & \mathbf{X}\mathbf{Z}^T \\ \mathbf{0} & \mu_u\mathbf{X}\mathbf{X}^T & \mathbf{0} \\ \mathbf{Z}\mathbf{X}^T & \mathbf{0} & \mathbf{Z}\mathbf{Z}^T \end{pmatrix}$$

and

$$B(\mathbf{X}, \mathbf{Z}) = \begin{pmatrix} \nu_0\mathbf{X}\mathbf{X}^T & \mathbf{0} & \nu_0\mathbf{X}\mathbf{Z}^T \\ 0 & \nu_2\mathbf{X}\mathbf{X}^T & 0 \\ \nu_0\mathbf{Z}\mathbf{X}^T & 0 & \nu_0\mathbf{Z}\mathbf{Z}^T \end{pmatrix}.$$

Denote the local likelihood estimate in (2.2) by $\tilde{\mathbf{a}}, \tilde{\mathbf{b}}$ and $\tilde{\boldsymbol{\beta}}$ and let

$$\hat{\boldsymbol{\beta}}^* = \sqrt{nh}\{(\tilde{\mathbf{a}} - \boldsymbol{\alpha}_0(u))^T, h(\tilde{\mathbf{b}} - \boldsymbol{\alpha}_0'(u))^T, (\tilde{\boldsymbol{\beta}} - \boldsymbol{\beta}_0)^T\}^T.$$

We then have the following asymptotic representation of $\widehat{\boldsymbol{\beta}}^*$.

LEMMA 2. *Under the regularity conditions given in Section* 5, *if* $h \to 0$ *and* $nh \to \infty$ *as* $n \to \infty$, *then* $\widehat{\boldsymbol{\beta}}^* = \mathbf{A}^{-1}\mathbf{W}_n + O_P\{h^2 + c_n \log^{1/2}(1/h)\}$ *holds uniformly in* $u \in \Omega$, *the support of* $U$, *where*

$$\mathbf{W}_n = \sqrt{h/n} \sum_{i=1}^{n} q_1(\bar{\alpha}_i, Y_i) \mathbf{X}_i^* K_h(U_i - u)$$

*and*

$$\mathbf{A} = f(u)E[\rho_2(\boldsymbol{\alpha}_0^T(U)\mathbf{X} + \mathbf{Z}^T\boldsymbol{\beta}_0)A(\mathbf{X}, \mathbf{Z})|U = u].$$

By some direct calculation, we then have the following mean and variance of $\mathbf{W}_n$:

$$E\mathbf{W}_n = \sqrt{nh}\frac{\mu_u}{2}\boldsymbol{\alpha}_0''^T(u)h^2 f(u)E[\rho_2\{\boldsymbol{\alpha}_0^T(U)\mathbf{X} + \mathbf{Z}^T\boldsymbol{\beta}_0\}(\mathbf{X}^T, \mathbf{0}, \mathbf{Z}^T)^T\mathbf{X}|U = u]$$
$$+ o(c_n^{-1}h^2)$$

and

$$\text{var}(\mathbf{W}_n) = f(u)E[\rho_2\{\boldsymbol{\alpha}_0^T(U)\mathbf{X} + \mathbf{Z}^T\boldsymbol{\beta}_0\}B(\mathbf{X}, \mathbf{Z})|U = u] + o(1).$$

Since $\mathbf{W}_n$ is a sum of independent and identically distributed random vectors, the asymptotic normality of $\tilde{\mathbf{a}}$, $\tilde{\mathbf{b}}$ and $\tilde{\boldsymbol{\beta}}$ can be established by using the central limit theorem and the Slutsky theorem. The next two theorems show that the estimate $\tilde{\boldsymbol{\beta}}$ can be improved by maximizing the penalized likelihood (2.3).

PROOF OF LEMMA 2. Throughout this proof, terms of the form $\widehat{G}(u) = O_P(a_n)$ always stand for $\sup_{u \in \Omega}|\widehat{G}(u)| = O_P(a_n)$. Let $c_n = (nh)^{-1/2}$. If $(\tilde{\mathbf{a}}, \tilde{\mathbf{b}}, \tilde{\boldsymbol{\beta}})^T$ maximizes (2.2), then $\widehat{\boldsymbol{\beta}}^*$ maximizes

$$\ell_n(\boldsymbol{\beta}^*) = h\sum_{i=1}^{n}[Q\{g^{-1}(c_n\boldsymbol{\beta}^{*T}\mathbf{X}_i^* + \bar{\alpha}_i), Y_i\} - Q\{g^{-1}(\bar{\alpha}_i), Y_i\}]K_h(U_i - u)$$

with respect to $\boldsymbol{\beta}^*$. The concavity of the function $\ell_n(\boldsymbol{\beta}^*)$ is ensured by condition (i). By a Taylor expansion of the function $Q\{g^{-1}(\cdot), Y_i\}$, we obtain that

$$(5.2) \qquad \ell_n(\boldsymbol{\beta}^*) = \mathbf{W}_n^T\boldsymbol{\beta}^* + \tfrac{1}{2}\boldsymbol{\beta}^{*T}\mathbf{A}_n\boldsymbol{\beta}^*\{1 + o_P(1)\},$$

where $\mathbf{A}_n = hc_n^2 \sum_{i=1}^{n} q_2(\bar{\alpha}_i, Y_i)\mathbf{X}_i^*\mathbf{X}_i^{*T}K_h(U_i - u)$. Furthermore, it can be shown that

$$(5.3) \qquad \mathbf{A}_n = -\mathbf{A} + o_P(1).$$

Therefore, by (5.2),

$$\ell_n(\boldsymbol{\beta}^*) = \mathbf{W}_n^T \boldsymbol{\beta}^* - \tfrac{1}{2} \boldsymbol{\beta}^{*T} \mathbf{A} \boldsymbol{\beta}^* + o_P(1). \tag{5.4}$$

Note that each element in $\mathbf{A}_n$ is a sum of i.i.d. random variables of kernel form and hence, by Lemma 1, converges uniformly to its corresponding element in $\mathbf{A}$. Consequently, expression (5.4) holds uniformly in $u \in \Omega$. By the Convexity Lemma (Pollard [19]), it also holds uniformly in $\boldsymbol{\beta}^* \in C$ and $u \in \Omega$ for any compact set $C$. Lemma A.1 of Carroll et al. [5] then yields

$$\sup_{u \in \Omega} |\widehat{\boldsymbol{\beta}}^* - \mathbf{A}^{-1} \mathbf{W}_n| \xrightarrow{P} 0. \tag{5.5}$$

Furthermore, from the definition of $\widehat{\boldsymbol{\beta}}^*$, we have that

$$\frac{\partial}{\partial \boldsymbol{\beta}^*} \ell_n(\boldsymbol{\beta}^*)\Big|_{\boldsymbol{\beta}^* = \widehat{\boldsymbol{\beta}}^*} = c_n h \sum_{i=1}^n q_1(\bar{\alpha}_i + c_n \widehat{\boldsymbol{\beta}}^{*T} \mathbf{X}_i^*, Y_i) \mathbf{X}_i^* \mathbf{X}_i^{*T} K_h(U_i - u) \widehat{\boldsymbol{\beta}}^* = 0.$$

By using (5.5) and a Taylor expansion, we have

$$\mathbf{W}_n + \mathbf{A}_n \widehat{\boldsymbol{\beta}}^* + \frac{c_n^3 h}{2} \sum_{i=1}^n q_3(\bar{\alpha}_i + \widehat{\zeta}_i, Y_i) \mathbf{X}_i^* \{\widehat{\boldsymbol{\beta}}^{*T} \mathbf{X}_i^*\}^2 K_h(U_i - u) = 0, \tag{5.6}$$

where $\widehat{\zeta}_i$ is between 0 and $c_n \widehat{\boldsymbol{\beta}}^{*T} \mathbf{X}_i^*$. The last term in the above expression is of order $O_P(c_n \|\widehat{\beta}^*\|^2)$. Since each element in $\mathbf{A}_n$ is of kernel form, we can deduce from Lemma 1 that $\mathbf{A}_n = E\mathbf{A}_n + O_P\{c_n \log^{1/2}(1/h)\} = -\mathbf{A} + O_P\{h^2 + c_n \log^{1/2}(1/h)\}$. Consequently, by (5.6), we obtain that

$$\mathbf{W}_n - \mathbf{A} \widehat{\boldsymbol{\beta}}^*[1 + O_P\{h^2 + c_n \log^{1/2}(1/h)\}] + O_P(c_n \|\widehat{\boldsymbol{\beta}}^*\|^2) = 0.$$

Hence,

$$\widehat{\boldsymbol{\beta}}^* = \mathbf{A}^{-1} \mathbf{W}_n + O_P\{h^2 + c_n \log^{1/2}(1/h)\}$$

holds uniformly for $u \in \Omega$. This completes the proof. $\quad\square$

PROOF OF THEOREM 1. Let $\gamma_n = n^{-1/2} + a_n$. It suffices to show that for any given $\zeta > 0$, there exists a large constant $C$ such that

$$P\left\{\sup_{\|\mathbf{v}\| = C} \mathcal{L}_P(\boldsymbol{\beta}_0 + \gamma_n \mathbf{v}) < \mathcal{L}_P(\boldsymbol{\beta}_0)\right\} \geq 1 - \zeta. \tag{5.7}$$

Define

$$D_{n,1} = \sum_{i=1}^n [Q\{g^{-1}(\hat{\boldsymbol{\alpha}}^T(U_i)\mathbf{X}_i + \mathbf{Z}_i^T(\boldsymbol{\beta}_0 + \gamma_n \mathbf{v})), Y_i\} $$
$$- Q\{g^{-1}(\hat{\boldsymbol{\alpha}}^T(U_i)\mathbf{X}_i + \mathbf{Z}_i^T \boldsymbol{\beta}_0), Y_i\}]$$

and $D_{n,2} = -n \sum_{j=1}^s \{p_{\lambda_n}(|\beta_{j0} + \gamma_n v_j|) - p_{\lambda_n}(|\beta_{j0}|)\}$, where $s$ is the number of components of $\boldsymbol{\beta}_{10}$. Note that $p_{\lambda_n}(0) = 0$ and $p_{\lambda_n}(|\beta|) \geq 0$ for all $\beta$. Thus,

$$\mathcal{L}_P(\boldsymbol{\beta}_0 + \gamma_n \mathbf{v}) - \mathcal{L}_P(\boldsymbol{\beta}_0) \leq D_{n,1} + D_{n,2}.$$

We first deal with $D_{n,1}$. Let $\widehat{m}_i = \widehat{\alpha}^T(U_i)\mathbf{X}_i + \mathbf{Z}_i^T \boldsymbol{\beta}_0$. Thus,

$$(5.8) \qquad D_{n,1} = \sum_{i=1}^n [Q\{g^{-1}(\widehat{m}_i + \gamma_n \mathbf{v}^T \mathbf{Z}_i), Y_i\} - Q\{g^{-1}(\widehat{m}_i), Y_i\}].$$

By means of a Taylor expansion, we obtain

$$(5.9) \qquad D_{n,1} = \sum_{i=1}^n q_1(\widehat{m}_i, Y_i)\gamma_n \mathbf{v}^T \mathbf{Z}_i - \frac{n}{2}\gamma_n^2 \mathbf{v}^T \mathbf{B}_n \mathbf{v},$$

where $\mathbf{B}_n = n^{-1} \sum_{i=1}^n \rho_2\{g^{-1}(\widehat{m}_i + \zeta_{ni})\}\mathbf{Z}_i \mathbf{Z}_i^T$, with $\zeta_{ni}$ between 0 and $\gamma_n \mathbf{v}^T \mathbf{Z}_i$, independent of $Y_i$. It can be shown that

$$(5.10) \quad \mathbf{B}_n = -E\rho_2\{\boldsymbol{\alpha}_0^T(U)\mathbf{X} + \mathbf{Z}^T\boldsymbol{\beta}_0\}\mathbf{Z}\mathbf{Z}^T + o_P(1) \equiv -\mathbf{B} + o_P(1).$$

Let $m_i = \boldsymbol{\alpha}_0^T(U_i)\mathbf{X}_i + \mathbf{Z}_i^T\boldsymbol{\beta}_0$. We have

$$n^{-1/2} \sum_{i=1}^n q_1(\widehat{m}_i, Y_i)\mathbf{Z}_i$$

$$= n^{-1/2} \sum_{i=1}^n q_1(m_i, Y_i)\mathbf{Z}_i$$

$$\qquad + n^{-1/2} \sum_{i=1}^n q_2(m_i, Y_i)[\{\widehat{\boldsymbol{\alpha}}(U_i) - \boldsymbol{\alpha}_0(U_i)\}^T \mathbf{X}_i]\mathbf{Z}_i$$

$$\qquad + O_P(n^{1/2}\|\widehat{\boldsymbol{\alpha}} - \boldsymbol{\alpha}_0\|_\infty^2).$$

By Lemma 2, the second term in the above expression can be expressed as

$$n^{-3/2} \sum_{i=1}^n q_2(m_i, Y_i)f(U_i)^{-1} \sum_{j=1}^n (\widetilde{W}_j^T \mathbf{X})K_h(U_j - U_i)\mathbf{Z}_i$$

$$+ O_P\{n^{1/2}c_n^2 \log^{1/2}(1/h)\}$$

$$\equiv T_{n1} + O_P\{n^{1/2}c_n^2 \log^{1/2}(1/h)\},$$

where $\widetilde{W}_j$ is the vector consisting of the first $p$ elements of $q_1(m_j, y_j)\Sigma^{-1}(u)$.

Define $\boldsymbol{\tau}_j = \boldsymbol{\tau}(\mathbf{X}_j, Y_j, \mathbf{Z}_j)$, consisting of the first $p$ elements of

$$q_1(m_j, Y_j)\boldsymbol{\Sigma}^{-1}(u)(\mathbf{X}_i^T, \mathbf{Z}_j^T)^T.$$

Using the definition of $\bar{\alpha}_j(U_i)$, we obtain $\bar{\alpha}_j(U_i) - m_j = O((U_j - U_i)^2)$ and therefore

$$T_{n1} = n^{-3/2} \sum_{i=1}^{n} \sum_{j=1}^{n} q_2(m_i, Y_i) f(U_i)^{-1} (\boldsymbol{\tau}_j^T \mathbf{X}_i) K_h(U_j - U_i) \mathbf{Z}_i + O_P(n^{1/2} h^2)$$

$$\equiv T_{n2} + O_P(n^{1/2} h^2).$$

It can be shown, by calculating the second moment, that

$$(5.11) \qquad\qquad T_{n2} - T_{n3} \xrightarrow{P} 0,$$

where $T_{n3} = -n^{-1/2} \sum_{j=1}^{n} \boldsymbol{\gamma}(U_j)$ with

$$\boldsymbol{\gamma}(u_j) = \sum_{k=1}^{p} \tau_{jk} E[\rho_2 \{\boldsymbol{\alpha}_0^T(u) \mathbf{X} + \mathbf{Z}^T \boldsymbol{\beta}_0\} X_k \mathbf{Z} | U = u_j].$$

Combining (5.8)–(5.11), we obtain that

$$(5.12) \qquad D_{n,1} = \gamma_n \mathbf{v} \sum_{i=1}^{n} \Omega(X_i, Y_i, \mathbf{Z}_i) - \tfrac{1}{2} \gamma_n^2 \mathbf{v}^T \mathbf{B} \mathbf{v} + o_P(1),$$

where $\Omega(U_i, Y_i, \mathbf{Z}_i) = q_1(m_i, Y_i) \mathbf{Z}_i - \boldsymbol{\gamma}(U_i)$. The orders of the first term and the second term are $O_P(n^{1/2} \gamma_n)$ and $O_P(n \gamma_n^2)$, respectively. We next deal with $D_{n,2}$. Note that $n^{-1} D_{n,2}$ is bounded by

$$\sqrt{s} \gamma_n a_n \|\mathbf{v}\| + \gamma_n^2 b_n \|\mathbf{v}\|^2 = C \gamma_n^2 (\sqrt{s} + b_n C),$$

by the Taylor expansion and the Cauchy–Schwarz inequality. As $b_n \to 0$, the second term on the right-hand side of (5.12) dominates $D_{n,2}$ as well as the first term on the right-hand side of (5.12), provided $C$ is taken to be sufficiently large. Hence, (5.7) holds for sufficiently large $C$. This completes the proof of the theorem. $\quad\square$

To prove Theorem 2, we need the following lemma.

LEMMA 3. *Under the conditions of Theorem* 2, *with probability tending to* 1, *for any given* $\boldsymbol{\beta}_1$ *satisfying* $\|\boldsymbol{\beta}_1 - \boldsymbol{\beta}_{10}\| = O_P(n^{-1/2})$ *and any constant* $C$, *we have*

$$\mathcal{L}_P \left\{ \begin{pmatrix} \boldsymbol{\beta}_1 \\ \mathbf{0} \end{pmatrix} \right\} = \max_{\|\boldsymbol{\beta}_2\| \leq C n^{-1/2}} \mathcal{L}_P \left\{ \begin{pmatrix} \boldsymbol{\beta}_1 \\ \boldsymbol{\beta}_2 \end{pmatrix} \right\}.$$

PROOF. We will show that, with probability tending to 1, as $n \to \infty$, for any $\boldsymbol{\beta}_1$ satisfying $\|\boldsymbol{\beta}_1 - \boldsymbol{\beta}_{10}\| = O_P(n^{-1/2})$ and $\boldsymbol{\beta}_2$ satisfying $\|\boldsymbol{\beta}_2\| \leq C n^{-1/2}$, $\partial \mathcal{L}_P(\boldsymbol{\beta}) / \partial \beta_j$ and $\beta_j$ have different signs for $\beta_j \in (-C n^{-1/2}, C n^{-1/2})$, for $j = s+1, \ldots, d$. Thus, the maximum is attained at $\boldsymbol{\beta}_2 = 0$. It follows by an argument similar to the proof of Theorem 1 that

$$\ell'_j(\boldsymbol{\beta}) \equiv \frac{\partial \ell(\tilde{\boldsymbol{\alpha}}, \boldsymbol{\beta})}{\partial \beta_j} = n \left\{ \frac{1}{n} \sum_{i=1}^n \Omega_j(\mathbf{X}_i, Y_i, \mathbf{Z}_i) - (\boldsymbol{\beta} - \boldsymbol{\beta}_0)^T B_j + o_P(n^{-1/2}) \right\},$$

where $\Omega_j(\mathbf{X}_i, Y_i, \mathbf{Z}_i)$ is the $j$th element of $\Omega(\mathbf{X}_i, Y_i, \mathbf{Z}_i)$ and $B_j$ is the $j$th column of $\mathbf{B}$. Note that $\|\boldsymbol{\beta} - \boldsymbol{\beta}_0\| = O_P(n^{-1/2})$ by the assumption. Thus, $n^{-1}\ell'_j(\boldsymbol{\beta})$ is of the order $O_P(n^{-1/2})$. Therefore, for $\beta_j \neq 0$ and $j = s+1, \ldots, d$,

$$\frac{\partial \mathcal{L}_P(\boldsymbol{\beta})}{\partial \beta_j} = \ell'_j(\boldsymbol{\beta}) - n p'_{\lambda_{jn}}(|\beta_j|) \operatorname{sgn}(\beta_j)$$

$$= -n\lambda_{jn} \left\{ \lambda_{jn}^{-1} p'_{\lambda_{jn}}(|\beta_j|) \operatorname{sgn}(\beta_j) + O_P\left(\frac{1}{\sqrt{n}\lambda_n}\right) \right\}.$$

Since $\liminf_{n \to \infty} \liminf_{\beta_j \to 0^+} \lambda_{jn}^{-1} p'_{\lambda_{jn}}(|\beta_j|) > 0$ and $\sqrt{n}\lambda_{jn} \to \infty$, the sign of the derivative is completely determined by that of $\beta_j$. This completes the proof. □

PROOF OF THEOREM 2. From Lemma 3, it follows that $\hat{\boldsymbol{\beta}}_2 = 0$. We next establish the asymptotic normality of $\hat{\boldsymbol{\beta}}_1$. Let $\hat{\boldsymbol{\theta}} = \sqrt{n}(\hat{\boldsymbol{\beta}}_1 - \boldsymbol{\beta}_{10})$, $\hat{m}_{i1} = \hat{\boldsymbol{\alpha}}^T(U_i)\mathbf{X}_i + \mathbf{Z}_{i1}^T\boldsymbol{\beta}_{10}$ and $m_{i1} = \boldsymbol{\alpha}_0^T(U_i)\mathbf{X}_i + \mathbf{Z}_{i1}^T\boldsymbol{\beta}_{10}$. Then, $\hat{\boldsymbol{\theta}}$ maximizes

$$\sum_{i=1}^n [Q\{g^{-1}(\hat{m}_{i1} + n^{-1/2}\mathbf{Z}_{i1}^T\boldsymbol{\theta}), Y_i\} - Q\{g^{-1}(\hat{m}_{i1}), Y_i\}] - n \sum_{j=1}^s p_{\lambda_n}(\hat{\beta}_{j1}).$$

(5.13)

We consider the first term, say $\ell_{n1}(\boldsymbol{\theta})$. By means of a Taylor expansion, we have

$$\ell_{n1}(\boldsymbol{\theta}) = n^{-1/2} \sum_{i=1}^n q_1(\hat{m}_{i1}, Y_i)\mathbf{Z}_{i1}^T\boldsymbol{\theta} + \frac{1}{2}\boldsymbol{\theta}^T \mathbf{B}_{n1}\boldsymbol{\theta},$$

where $\mathbf{B}_{n1} = \frac{1}{n} \sum_{i=1}^n \rho_2\{g^{-1}(\hat{m}_{i1} + \zeta_{ni})\}\mathbf{Z}_{i1}\mathbf{Z}_{i1}^T$, with $\zeta_{ni}$ between 0 and $n^{-1/2} \times \mathbf{Z}_{i1}^T\boldsymbol{\theta}$, independent of $Y_i$. It can be shown that

(5.14)   $\mathbf{B}_{n1} = -E\rho_2\{\boldsymbol{\alpha}_0^T(U)\mathbf{X} + \mathbf{Z}_1^T\boldsymbol{\beta}_{10}\}\mathbf{Z}_1\mathbf{Z}_1^T + o_P(1) = -\mathbf{B}_1 + o_P(1).$

A similar proof for (5.12) yields that

$$\ell_{n1}(\boldsymbol{\theta}) = n^{-1/2} \sum_{i=1}^n \hat{\boldsymbol{\theta}} \Omega_1(U_i, Y_i, \mathbf{Z}_{i1}) - \frac{1}{2}\boldsymbol{\theta}^T \mathbf{B}_1 \boldsymbol{\theta} + o_P(1),$$

where $\Omega_1(U_i, Y_i, \mathbf{Z}_{i1}) = q_1(m_{i1}, Y_i)\mathbf{Z}_{i1} - \Gamma_1(U_i)$. By the Convexity Lemma (Pollard [19]), we have that

$$(\mathbf{B}_1 + \Sigma_\lambda)\hat{\boldsymbol{\theta}} + n^{1/2}\mathbf{b}_n = n^{-1/2} \sum_{i=1}^n \Omega_1(U_i, Y_i, \mathbf{Z}_{i1}) + o_P(1).$$

The conclusion follows as claimed. □

PROOF OF THEOREM 3. Decompose $\mathcal{R}(H_1) - \mathcal{R}(H_0) = I_{n,1} + I_{n,2} + I_{n,3}$, where

$$I_{n,1} = \sum_{i=1}^{n} [Q\{g^{-1}(\widehat{\boldsymbol{\alpha}}^T(U_i)\mathbf{X}_i + \mathbf{Z}_i^T\widehat{\boldsymbol{\beta}}), Y_i\} - Q\{g^{-1}(\widehat{\boldsymbol{\alpha}}^T(U_i)\mathbf{X}_i + \mathbf{Z}_i^T\boldsymbol{\beta}_0), Y_i\}],$$

$$I_{n,2} = -\sum_{i=1}^{n} [Q\{g^{-1}(\mathbf{Z}_i^T\bar{\boldsymbol{\beta}}), Y_i\} - Q\{g^{-1}(\mathbf{Z}_i^T\boldsymbol{\beta}_0), Y_i\}],$$

$$I_{n,3} = \sum_{i=1}^{n} [Q\{g^{-1}(\widehat{\boldsymbol{\alpha}}^T(U_i)\mathbf{X}_i + \mathbf{Z}_i^T\boldsymbol{\beta}_0), Y_i\} - Q\{g^{-1}(\mathbf{Z}_i^T\boldsymbol{\beta}_0), Y_i\}].$$

Using Theorem 10 of Fan, Zhang and Zhang [10], under $H_0$, we have

$$r_K I_{n,3} \sim \chi^2_{\mathrm{df}_n},$$

where $\mathrm{df}_n \to \infty$ as $n \to \infty$. It suffices to show that $I_{n,1} = o_P(I_{n,3})$ and $I_{n,2} = o_P(I_{n,3})$.

A direct calculation yields that

$$I_{n,1} = \sum_{i=1}^{n} q_1\{\mathbf{X}_i^T\widehat{\boldsymbol{\alpha}}(U_i) + \mathbf{Z}_i^T\boldsymbol{\beta}_0, Y_i\}\mathbf{Z}_i^T(\widehat{\boldsymbol{\beta}} - \boldsymbol{\beta}_0)$$

$$- \tfrac{1}{2}(\widehat{\boldsymbol{\beta}} - \boldsymbol{\beta}_0)^T \sum_{i=1}^{n} \mathbf{Z}_i\mathbf{Z}_i^T q_2\{g^{-1}(\widehat{\alpha}(U_i)\mathbf{X}_i + \mathbf{Z}_i^T\boldsymbol{\beta}_0)\}(\widehat{\boldsymbol{\beta}} - \boldsymbol{\beta}_0) + o_p(1).$$

Using techniques related to those used in the proof of Theorem 2, we obtain

$$\frac{1}{n}\sum_{i=1}^{n} \mathbf{Z}_i\mathbf{Z}_i^T q_2[g^{-1}\{\widehat{\alpha}(U_i)\mathbf{X}_i + \mathbf{Z}_i^T\boldsymbol{\beta}_0\}] = \mathbf{B} + o_P(1),$$

$$\sum_{i=1}^{n} q_1\{\widehat{\boldsymbol{\alpha}}(U_i)\mathbf{X}_i + \mathbf{Z}_i^T\boldsymbol{\beta}_0, Y_i\}\mathbf{Z}_i = n\mathbf{B}(\hat{\beta} - \boldsymbol{\beta}_0) + o_P(1).$$

Thus,

$$2I_{n,1} = (\hat{\boldsymbol{\beta}} - \boldsymbol{\beta}_0)^T \mathbf{B}(\hat{\boldsymbol{\beta}} - \boldsymbol{\beta}_0) + o_p(1) \xrightarrow{D} \chi^2_d.$$

Under $H_0$, $-2I_{n,2}$ equals a likelihood ratio test statistic for $H_0^*: \boldsymbol{\beta} = \boldsymbol{\beta}_0$ versus $H_1^*: \boldsymbol{\beta} \neq \boldsymbol{\beta}_0$. Thus, under $H_0$, $-2I_{n,2} \to \chi^2_d$. Thus, $I_{n,1} = o_P(I_{n,3})$ and $I_{n,2} = o_P(I_{n,3})$. This completes the proof. □

**Acknowledgments.** The authors thank the Co-Editor and the referees for constructive comments that substantially improved an earlier version of this paper.

DEPARTMENT OF STATISTICS
  AND THE METHODOLOGY CENTER
PENNSYLVANIA STATE UNIVERSITY
UNIVERSITY PARK, PENNSYLVANIA 16802-2111
USA
E-MAIL: rli@stat.psu.edu

DEPARTMENT OF BIOSTATISTICS
  AND COMPUTATIONAL BIOLOGY
UNIVERSITY OF ROCHESTER
ROCHESTER, NEW YORK 14642
USA
E-MAIL: hliang@bst.rochester.edu